\documentclass{elsart}
\usepackage{graphicx}
\usepackage{graphics}
\usepackage[dvips]{epsfig}
\journal{Statistics and Probability Letters}

\newcommand{\bxi}{{\mbox{\boldmath $\xi$}}}

\newcommand{\bU}{{\bf U}}

\newcommand{\bR}{{\bf R}}

\newcommand{\N}{I \! \! {N}}
\newcommand{\C}{I \! \! \! \! {C}}

\newcommand{\umu}{{\underline \mu}}

\newcommand{\ux}{{\underline x}}
\newcommand{\uy}{{\underline y}}

\newcommand{\p}{{p^*}}

\newcommand{\bb}{\begin{eqnarray}}
\newcommand{\be}{\end{eqnarray}}

\newtheorem{theorem}{Theorem}

\newtheorem{lemma}{Lemma}

\begin{document}
\begin{frontmatter}

\title{On the condensed density of the zeros of the Cauchy transform of a complex atomic random measure with Gaussian moments}
\author{P. Barone }
\address{ Istituto per le Applicazioni del Calcolo ''M. Picone'',
C.N.R.,\\
Via dei Taurini 19, 00185 Rome, Italy \\
e-mail: piero.barone@gmail.com, p.barone@iac.cnr.it  \\
fax: 39-6-4404306}

%\maketitle

\begin{abstract}

%\section*{Abstract}

An atomic random complex measure defined on the unit disk with Normally distributed moments is considered. An approximation to the distribution of the zeros of its Cauchy transform is computed.
Implications of this result for solving several moments problems are discussed.

\end{abstract}


\end{frontmatter}
%\newpage

\section*{Introduction}

Let us consider the complex random measure
$${\bf S}(z)=\sum_{j=1}^{p}{\bf c}_j\delta(z-\mbox{\boldmath $\xi$}_j),\;z\in  D,\;\;Prob[\mbox{\boldmath $\xi$}_j\in D]=1$$
where $D$ denotes the unit disk and bold characters denote random quantities, and assume that all finite sets of its complex moments
\begin{eqnarray}{\bf d}_k=\int_{ D} z^k d{\bf S}(z)=\sum_{j=1}^{p}{\bf c}_j\mbox{\boldmath $\xi$}_j^k,\;\;k\in \N.\label{eq1}\end{eqnarray}
have a joint complex Gaussian distribution
$$g(\ux)=\frac{1}{(\pi \sigma^2)^n}e^{-\frac{1}{\sigma^2}\sum_{k=0}^{n-1}|x_k-\mu_k|^2}$$
where $\mu_k=E[{\bf d}_k]\in\C.$
The Cauchy transform of ${\bf S}$ is defined as
$${\bf f}(z)=\int_{D} \frac{d{\bf S}(\zeta)}{z-\zeta}=\sum_{j=1}^{p}\frac{{\bf c}_j}{z-\mbox{\boldmath $\xi$}_j}.$$
We are looking for the distribution of the zeros of ${\bf f}(z)$.

Some motivations are provided in the following.  We notice that, in the deterministic case, $f(z)$ is equal to the derivative of the logarithmic potential
$$L(z)=\log\prod_{j=1}^{p}(z-\xi_j)^{ c_j}.$$
Therefore  the zeros of the Cauchy transform are the stationary points of $L(z)$ i.e. the location of the equilibrium points in a field of force due to complex masses $c_j$ at the points $\xi_j$ acting according to the inverse distance law in the plane. These locations were described in \cite[Th.(8,2)]{ma} as a generalization of Lucas' theorem \cite[Th.(6,1)]{ma}.

It turns out (see e.g. \cite{barja2,bardsp} that many difficult inverse problems in a stochastic framework, can be reduced to the estimation of a measure ${\bf S}(z)$ from its moments ${\bf d}_k$. This is equivalent to make inference on $p,\;{\bf c}_j,\;\mbox{\boldmath $\xi$}_j$ from ${\bf d}_k$. As   ${\bf S}(z)$ depends in a highly non linear way on $p$ and $\mbox{\boldmath $\xi$}_j$, these are the most critical quantities to estimate. In \cite{barja2} an approach to cope with this problem was proposed which is based on the estimation of the condensed density of the $\mbox{\boldmath $\xi$}_j$, i.e. the poles of ${\bf f}(z)$, which is defined as
$$h(z)= \frac{1}{p}E\left[\sum_{j=1}^{p}\delta(z-\bxi_j)\right]$$ or,
equivalently, for all Borel sets $A\subset\C$ $$\int_A
h(z)dz=\frac{1}{p}\sum_{j=1}^p Prob(\bxi_j\in A).$$ It can be proved
that (see e.g. \cite{barja})
$$h(z)=\frac{1}{4\pi}\Delta u(z) $$ where
$\Delta$ denotes the Laplacian operator  with respect to $x,y$ if
$z=x+iy$ and
$$u(z)=\frac{1}{p}E\left\{\log(|\prod_{j=1}^p(z-\mbox{\boldmath $\xi$}_j)|^2)\right\}$$ is the
corresponding logarithmic potential.

However the zeros of ${\bf f}(z)$, after the logarithmic potential interpretation, convey information about the $\mbox{\boldmath $\xi$}_j$. In the simplest deterministic case of positive $c_j$ and $p=2$, the only  zero of $f(z)$ is the barycenter of  $\xi_1,\xi_2$ to which the weights $\frac{c_2}{c_1+c_2},\frac{c_1}{c_1+c_2}$ are associated. In the general case a qualitative picture of the zeros location w.r. to the $\xi_j$ is the following: to each couple of $\xi_h,\xi_k$ it corresponds a zero in a strip connecting $\xi_h,\xi_k$. The strip shape and the position of the zero in the strip is determined mainly by  $c_h,c_k$: the smaller  $|c_h|$ is w.r. to $|c_k|$ the closest  the zero is to $\xi_h$. Therefore knowledge of the zeros provide constraints to $\xi_j$ and $c_j$ useful for estimating them (see fig.1).

In the following a method to estimate the condensed density of the zeros of  ${\bf f}(z)$ starting from simple statistics which can be computed from a sample of the moments ${\bf d}_k$ is given. The idea is to interpret the zeros of ${\bf f}(z)$ as the poles of the Cauchy transform of a random measure ${\bf \tilde{S}}(z)$ associated to ${\bf S}(z)$ through a deterministic one-to-one transformation and to apply the same approximation technique used for estimating the condensed density $h(z)$ of the poles $\mbox{\boldmath $\xi$}_j$ in \cite{distrf,barspl}.

\section{Pade' approximants}

In \cite{distrf}, in order to compute the condensed density of the poles of ${\bf f}(z)$ assuming that the joint distribution of a finite number of moments ${\bf d}_0,\dots,{\bf d}_{2p-1}$ is Gaussian, a  pencil of square random Hankel matrices $$\bU(z)=\bU_1-z\bU_0$$ was built from ${\bf d}_0,\dots,{\bf d}_{2p-1}$ such that $\mbox{\boldmath $\xi$}_j$ are its generalized eigenvalues:
$$\bU_0=\left[\begin{array}{llll}
{\bf d}_0 & {\bf d}_{1} &\dots &{\bf d}_{p-1} \\
{\bf d}_{1} & {\bf d}_{2} &\dots &{\bf d}_{p} \\
. & . &\dots &. \\
{\bf d}_{p-1} & {\bf d}_{p} &\dots &{\bf d}_{2p-2}
  \end{array}\right],\;
\bU_1=\left[\begin{array}{llll}
{\bf d}_1 & {\bf d}_{2} &\dots &{\bf d}_{p} \\
{\bf d}_{2} & {\bf d}_{3} &\dots &{\bf d}_{p+1} \\
. & . &\dots &. \\
{\bf d}_{p} & {\bf d}_{p+1} &\dots &{\bf d}_{2p-1}\end{array}\right]$$
This property follows from equation (\ref{eq1})  (see e.g. \cite[Sec.7.2]{hen2}, \cite{barja2}).
But then
\begin{eqnarray}u(z)&=&\frac{1}{p}E\left\{\log(|\prod_{j=1}^p(z-\mbox{\boldmath $\xi$}_j)|^2)\right\}=
E\left\{\log(|\det[\bU_0^{-1}\bU_1-zI_p]|^2)\right\}\nonumber\\&=&E\left\{\log(|\det[\bU(z)]|^2|\det[\bU_0^{-1}]|^2)\right\}\nonumber.\end{eqnarray}
By noticing that  $\det[\bU_0^{-1}]|^2$ does not depend on $z$, and by taking the $QR$ factorization of $\bU(z)$ given by Gram-Schmidt algorithm we got
  $$h(z)=\frac{1}{4\pi}\Delta u(z)= \frac{1}{4\pi p}\Delta
\sum_{k=1}^pE\left\{\log(\bR_{kk}(z))^2\right\} .$$ The distribution of
$\bR_{kk}^2,k=1,\dots,p$ was then expanded in an uniformly convergent  series of generalized Laguerre functions and an explicit expression for
$E[\log(\bR_{kk})^2]$ and then for  $h(z)$ was found. By truncating the series a computable approximation of $h(z)$ was derived as a function of the moments of $\bR_{kk}^2$.

 In order to apply the methods developed in
\cite{distrf} we look for a new moment sequence ${\bf \tilde{d}}_k$
$${\bf \tilde{d}}_k=\sum_{j=1}^{p-1}{\bf \tilde{c}}_j\mbox{\boldmath $\zeta$}_j^k,\;\;k\in \N $$
and the associated measure
$${\bf \tilde{S}}(z)=\sum_{j=1}^{p-1}{\bf \tilde{c}}_j\delta(z-\mbox{\boldmath $\zeta$}_j)$$ such that
 $\mbox{\boldmath $\zeta$}_j$ are the zeros of ${\bf f}(z)$. The following theorem holds:
\begin{theorem}
 Let us define the first $k$ elements $\underline{{\bf \tilde{d}}}$ of the sequence ${\bf \tilde{d}}_k$ for all $k$  by \begin{eqnarray} \underline{{\bf \tilde{d}}}={\bf T}^{-1}(\underline{{\bf d}}) \underline{e_1} \label{eq3}\end{eqnarray}
  where
$${\bf T}={\bf T}(\underline{{\bf d}})=\left[\begin{array}{llll}
{\bf d}_0 & 0 &\dots &0 \\
{\bf d}_{1} & {\bf d}_{0} &\dots &0 \\
... \\
{\bf d}_{k-1} & {\bf d}_{k-2} &\dots &{\bf d}_{0}
  \end{array}\right],\;\;\;\underline{{\bf d}}=\left[\begin{array}{llll}
{\bf d}_0  \\
{\bf d}_{1}  \\
... \\
{\bf d}_{k-1}
  \end{array}\right],\;\;\;\underline{e_1}=\left[\begin{array}{llll}
1  \\
0  \\
... \\
0
  \end{array}\right].$$ Then
  the zeros of ${\bf f}(z)$ are the generalized eigenvalues of
 the $(p-1)\times (p-1)$ pencil ${\bf \tilde{\bU}}(z)={\bf \tilde{\bU}}_1-z{\bf \tilde{\bU}}_0$ where
$${\bf \tilde{\bU}}_0=\left[\begin{array}{llll}
{\bf \tilde{d}}_2 & {\bf \tilde{d}}_{3} &\dots &{\bf \tilde{d}}_{p} \\
{\bf \tilde{d}}_{3} & {\bf \tilde{d}}_{4} &\dots &{\bf \tilde{d}}_{p+1} \\
. & . &\dots &. \\
{\bf \tilde{d}}_{p} & {\bf \tilde{d}}_{p+1} &\dots &{\bf \tilde{d}}_{2p-2}
  \end{array}\right],\;
{\bf \tilde{\bU}}_1=\left[\begin{array}{llll}
{\bf \tilde{d}}_3 & {\bf \tilde{d}}_{4} &\dots &{\bf \tilde{d}}_{p+1} \\
{\bf \tilde{d}}_{4} & {\bf \tilde{d}}_{5} &\dots &{\bf \tilde{d}}_{p+2} \\
. & . &\dots &. \\
{\bf \tilde{d}}_{p+1} & {\bf \tilde{d}}_{p+2} &\dots &{\bf \tilde{d}}_{2p-1}\end{array}\right].$$ Moreover
  \begin{eqnarray}{\bf \tilde{d}}_k=\sum_{j=1}^{p-1}{\bf \tilde{c}}_j\mbox{\boldmath $\zeta$}_j^k,\;\;k\in \N \label{eq2}.\end{eqnarray}
\end{theorem}
\noindent\underline{proof.}
Let us consider the formal random power series
$${\bf F}(z)=\sum_{k=0}^{\infty}{\bf d}_kz^{-k}=\sum_{j=1}^{p}
{\bf c}_j\sum_{k=0}^{\infty}(\mbox{\boldmath $\xi$}_j/z)^k= z {\bf f}(z),\;\;|z|>1$$ which can be extended to $D$ by analytic continuation.
The Cauchy transform ${\bf f}(z)$ can then be seen as the  Pade' approximant of ${\bf F}(z)$ of orders $p-1,p$, denoted by $[p-1,p]_{\bf F}(z)$, times $z^{-1}$.
If $[p-1,p]_{\bf F}(z)=\frac{ {\bf Q}_{p-1}(z)}{{\bf P}_p(z)}$ where ${\bf P}_p(z), {\bf Q}_{p-1}(z)$ are polynomial of degree $p,p-1$ respectively, the polynomials ${\bf P}_k(z)$ are orthogonal w.r. to the sequence  $\{{\bf d}_0,{\bf d}_1,\dots\}$ \cite[Sec. 2.1]{bre} and therefore the zeros of ${\bf P}_k(z)$ are the generalized eigenvalues of the pencil $\bU(z)=\bU_1-z\bU_0$.

Let us
define the sequence  ${\bf \tilde{d}}_k$ by the identity
\begin{eqnarray}{\bf F}(z)^{-1}=\sum_{k=0}^{\infty}{\bf \tilde{d}}_kz^{-k}\label{fm1}\end{eqnarray} which implies that $\forall k\in \N$
$${\bf T}\underline{{\bf \tilde{d}}}= \underline{e_1} $$
The random triangular matrix ${\bf T}$ is a.s. invertible as ${\bf d}_0$ is a.s. different from zero. Hence ${\bf \tilde{d}}_k$ is well defined.
From \cite[Cor. 2.7]{bre} it follows that the polynomials ${\bf Q}_k(z)$ are orthogonal w.r. to the sequence  $\{{\bf \tilde{d}}_2,{\bf \tilde{d}}_3,\dots\}$. Therefore the roots of ${\bf Q}_{p-1}(z)$, i.e. the zeros of ${\bf f}(z)$, are the generalized eigenvalues of ${\bf \tilde{\bU}}(z).$

Moreover we  notice that the Pade' approximant $\mbox{$[p,p-1]_{{\bf F}^{-1}}(z)$}$ to the reciprocal power series $z{\bf F}^{-1}(z)$ is $\frac{1}{[p-1,p]_{\bf F}(z)}$. Therefore the zeros of $z {\bf f}(z)=[p-1,p]_{\bf F}(z)$ are the poles of $[p,p-1]_{{\bf F}^{-1}}(z)=\frac{{\bf P}_p(z)}{ {\bf Q}_{p-1}(z)}$.  Because  numerator has degree greater than  denominator we can divide  ${\bf P}_p(z)$ by ${\bf Q}_{p-1}(z)$ getting
$$ [p,p-1]_{{\bf F}^{-1}}(z)=\frac{{\bf P}_p(z)}{ {\bf Q}_{p-1}(z)}=\sum_{j=1}^{p-1}\frac{{\bf \tilde{c}}_j}{z-\mbox{\boldmath $\zeta$}_j}+{\bf D}_1(z)$$
where ${\bf D}_1(z)$ is a polynomial of degree one. But then $$[p-2,p-1]_{{\bf F}^{-1}}(z)=\sum_{j=1}^{p-1}\frac{{\bf \tilde{c}}_j}{z-\mbox{\boldmath $\zeta$}_j}$$ and from eq.\ref{fm1}, we get ${\bf \tilde{d}}_k=\sum_{j=1}^{p-1}{\bf \tilde{c}}_j\mbox{\boldmath $\zeta$}_j^k.\;\;\;\;\Box$

\section{Approximate distribution of the moments of the associated measure}

In order to derive the joint density of ${\bf \tilde{d}}_0,\dots,{\bf \tilde{d}}_{n-1}$ we need the following lemma
\begin{lemma}
The transformation $\Phi: \underline{{\bf d}}\rightarrow \underline{{\bf \tilde{d}}}$ defined by
\begin{eqnarray} \underline{{\bf \tilde{d}}}={\bf T}^{-1}(\underline{{\bf d}}) \underline{e_1} \label{eq4}\end{eqnarray}
is an involution i.e. $\Phi(\Phi({\bf d}))={\bf d}$.
\label{lem1}
\end{lemma}
\noindent\underline{proof.}
The thesis follows by the Toeplitz structure of the matrix ${\bf T}.\;\;\;\Box$
\begin{theorem}
The joint density of ${\bf \tilde{d}}_0,\dots,{\bf \tilde{d}}_{n-1}$
is given by
$$\tilde{g}(\uy)=\frac{1}{(\pi \sigma^2y_0^2)^n }e^{-\frac{1}{\sigma^2}\sum_{k=0}^{n-1}|[\Phi^{-1}(\uy)]_k-\mu_k|^2}$$
where
$$[\Phi^{-1}(\uy)]_k=\frac{(-1)^{k(k+1)/2}H_k^{(-k+2)}}{y_0^{k+1}}$$
and
$$H_k^{(-k+2)}=\det\left[\begin{array}{lllll}
y_1 & y_2 &\dots & y_{n-1}&y_n \\
y_0 & y_1 &\dots & y_{n-2}&y_{n-1} \\
. & . &\dots &. &.\\
0 & 0 &\dots& y_0 &y_1\end{array}\right].$$
\end{theorem}
\noindent\underline{proof.}
If $n=2p$ the joint density of ${\bf d}_0,\dots,{\bf d}_{n-1}$ is
$$g(\ux)=\frac{1}{(\pi \sigma^2)^n}e^{-\frac{1}{\sigma^2}\sum_{k=0}^{n-1}|x_k-\mu_k|^2}.$$
By Lemma \ref{lem1}, if $\uy=\Phi(\ux),\;\ux,\uy\in\C^n$, the complex Jacobian of the transformation $\Phi$
is
$$J_c=\frac{(-1)^n}{{y_0}^{2n}}=\frac{1}{{y_0}^{2n}}$$
because $n$ is even.
But then, the joint density of ${\bf \tilde{d}}_0,\dots,{\bf \tilde{d}}_{n-1}$ is
$$\tilde{g}(\uy)=\frac{1}{(\pi \sigma^2y_0^2)^n }e^{-\frac{1}{\sigma^2}\sum_{k=0}^{n-1}|[\Phi^{-1}(\uy)]_k-\mu_k|^2}.$$
The rest of the thesis follows by \cite[pg.106]{bre}. $\;\;\;\Box$
\begin{theorem}
In the limit for $\sigma\rightarrow 0$, $\tilde{g}(\uy)$ tends to a Dirac measure  centered in $\Phi(\umu)$.
\end{theorem}
\noindent\underline{proof.}
By noticing that $\tilde{g}(\uy)$ can be written as
$$\tilde{g}(\uy)=\int\delta(\uy-\Phi(\ux))g(\ux)d\ux$$
in the limit for $\sigma\rightarrow 0$ we get, for every continuous function $\psi(\uy)$ with compact support
$$\lim_{\sigma\rightarrow 0}\int\tilde{g}(\uy)\psi(\uy)d\uy=\lim_{\sigma\rightarrow 0}\int\int\delta(\uy-\Phi(\ux))g(\ux)\,\psi(\uy) d\ux d\uy=$$
$$\lim_{\sigma\rightarrow 0}\int g(\ux)\left(\int\delta(\uy-\Phi(\ux))\,\psi(\uy)d\uy\right)d\ux=\lim_{\sigma\rightarrow 0}\int g(\ux)\psi(\Phi(\ux))d\ux=$$
$$\int \delta(\ux-\umu)\psi(\Phi(\ux))d\ux=\psi(\Phi(\umu)).\;\;\;\Box$$
As a consequence of this theorem we can approximate $\tilde{g}(\uy)$, in the limit for $\sigma\rightarrow 0$, by a Gaussian  centered in $\Phi(\umu)$ with height
$$\tilde{g}(\Phi(\umu))=\frac{1}{(\pi \sigma^2[\Phi(\umu)]_0^2)^n }=\frac{1}{(\pi \frac{\sigma^2}{\umu_0^2})^n }$$
i.e.
$$\tilde{g}(\uy)\approx \frac{\umu_0^{2n}}{(\pi \sigma^2)^n}e^{-\frac{\umu_0^2}{\sigma^2}\sum_{k=0}^{n-1}|y_k-[\Phi(\umu)]_k|^2}.$$

\section{Approximate  condensed density of the zeros}

In view of the result of the previous section, we can use the theory developed in \cite{distrf,barspl} to find an approximation of the condensed density of the zeros of ${\bf f}(z)$.
More precisely  starting from new data ${\bf \tilde{d}}_0,\dots,{\bf \tilde{d}}_{n-1}$ that are approximately distributed as a Gaussian, we can compute the QR factorization of the pencil ${\bf \tilde{\bU}}(z).$ The squared diagonal elements $\bR_{kk}^2,k=1,\dots,p-1$ are then (conditionally) distributed as quadratic forms in Normal variables and then their distribution can be approximated by a random variable whose density $f_k(y)$ admits a uniformly convergent Laguerre series ( \cite[Th.1]{distrf})
 \begin{eqnarray*} f_k(y)=b^{(k)}_0\frac{y^{\alpha_k-1}e^{-y/\beta_k}}{\beta_k^{\alpha_k}\Gamma(\alpha_k)}+
\frac{y^{\alpha_k-1}e^{-y/\beta_k}}{\beta_k^{\alpha_k}\Gamma(\alpha_k)}\sum_{m=1}^\infty b^{(k)}_m L_m(y/\tau_k,\alpha_k).\end{eqnarray*}
 But then $E[\log(\bR_{kk})^2]$ can be approximated by taking the expectation of each term of the series which can be computed in closed form giving
$$E\left\{\log|\bR_{kk}(z)|^2\right\}\approx
 b^{(k)}_0(z)[\log\beta_k(z)+\Psi(\alpha_k(z))]+$$
$$\sum_{m=1}^\infty b^{(k)}_m(z)\sum_{h=0}^m c_{hm}\frac{\Gamma (\alpha_k(z)+h)} {\Gamma (\alpha_k(z))}\left(\frac{\beta_k(z)}{\tau_k(z)}\right)^h \left[\log \beta_k(z)+\Psi(\alpha_k(z)+h)\right] $$
where $c_{hm}$ are the coefficient of the Laguerre polynomials, $b^{(k)}_m(z)$ are the coefficient of the Laguerre expansion, $\alpha_k(z),\beta_k(z)$ are the parameters of the Gamma density which is the first term of the series for $f_k(y)$, and $\Psi(\cdot)$ is the logarithmic derivative of the Gamma function.

The estimate of the condensed density of the zeros of ${\bf f}(z)$ is then obtained by truncating the series after the first term,  replacing $\alpha_k(z)$ by convenient estimates, and $\beta_k(z)$ by   $\beta>0$ which turns out to be a smoothing parameter (\cite[Sec.3]{distrf}):
$$\hat{h}(z)\propto\sum_{k=1}^{p-1}\hat{\Delta}
\left(\Psi\left[\left(
\frac{\hat{R}_{kk}^2(z)}{\sigma^2\beta}+1\right)\right]\right)$$
where $\hat{\Delta}$ is the discrete Laplacian, and $\hat{R}_{kk}^2(z)$ is obtained by the QR factorization of a realization of the pencil ${\bf \tilde{\bU}}(z)$.

\section{Simulation results}

To appreciate the goodness of the approximation to the density of $\bR_{kk}^2$ provided by the truncated Laguerre expansion, for $n=74,\;N=10^6$  independent
  realizations $d_k^{(r)}, k=1,\dots,n,\;\;r=1,\dots,N$ of the r.v.  ${\bf d}_k$ were
  generated by the equation
  $$d_k^{(r)}=s_k+\epsilon_j^{(r)},\;\;s_k=\sum_{j=1}^\p c_j\xi_j^k$$ where $\epsilon_j^{(r)}$ are independent realizations of zero mean Gaussian variables with standard deviation $\sigma=0.2$,  $\p=5$ and
$$\underline{\xi}=\left[ e^{-0.1-i 2\pi  0.3},e^{-0.05-i 2\pi
0.28},e^{-0.0001+i 2\pi 0.2},e^{-0.0001+i 2\pi  0.21},e^{-0.3-i 2\pi
0.35}\right]$$ $$ \underline{c}=\left[ 6,3,1,1,20\right
].$$
The matrix $T^{(r)}$ based on $d_k^{(r)}$ was formed and the linear system (\ref{eq3}) was solved for the new moments
$\tilde{d}^{(r)}_k,k=1,\dots,n,\;\;r=1,\dots,N$. The $(p-1)\times (p-1)$ matrices $\tilde{U}^{(r)}_0,\;\tilde{U}^{(r)}_1$
based on $\tilde{d}_k^{(r)}$ were computed. The
 matrix $\tilde{U}_1^{(r)}-z\tilde{U}_0^{(r)}$ with $z=\cos(1)+i0.8$ was formed, its $QR$ decomposition  and the first $10$ empirical moments $\hat{\gamma}_j$ were computed. Estimates of the first $10$ coefficients of the Laguerre expansion were then computed by (\cite{balak})
 $$\hat{\alpha}_k= \frac{\hat{\gamma}_1^2}{\hat{\gamma}_2-\hat{\gamma}_1^2} ,\;\; \hat{\beta}_k=\frac{\hat{\gamma}_2-\hat{\gamma}_1^2}{\hat{\gamma}_1}  $$
 $$\hat{b}^{(k)}_h= (-1)^h\Gamma(\hat{\alpha}_k)\sum_{j=0}^h(-1)^j{h \choose j}\frac{\hat{\gamma}_{h-j}}{\Gamma(\hat{\alpha}_{k+h-j})},\;\;\hat{\gamma}_{0}=1,\;\;h=1,\dots,10$$
  The one term and ten terms approximations of the density were then computed and compared with the empirical density of $\bR_{kk}^2$ for $k=1,\dots,p$. The results are given in fig.\ref{fig0}. In the
top left part the real part of one realization of the transformed data $\tilde{d}_k$ and  transformed signal $\tilde{s}_k$ are plotted. In the
 top right part the   $L_2$ norm of the difference between the empirical density  of  $\bR_{kk}^2,k=1,\dots,p$  computed by MonteCarlo simulation and its approximation obtained by truncating the series expansion of the density after the first term and after the first $10$ terms is given.
 In the bottom left part the  density of $\bR_{kk}^2,\;k=35,$ approximated by the first term of its series expansion and the empirical density are plotted.
  In the bottom right part the density of $\bR_{kk}^2,\;k=35,$ approximated by the first $10$ terms of its series expansion and the empirical density are plotted.
  It can be noticed that the first order approximation is quite good even if it becomes worse for large $k$.

To appreciate the advantage of the closed form estimate $\hat{h}(z)$
with respect to an estimate of the condensed density obtained by
MonteCarlo simulation an
 experiment was performed. $N=100$ independent
  realizations of the  r.v. generated above were considered.

An estimate of $h(z)$ was computed on a
square lattice of dimension $m=100$ by
$$\hat{h}(z)\propto
\sum_{r=1}^N\sum_{k=1}^p\hat{\Delta}
\left\{\Psi\left[\left(\frac{R_{kk}^{(r)}(z)^2}{\sigma^2\beta}+1\right)\right]\right\}
$$ where $R^{(r)}(z)$
is obtained by the QR factorization of the matrix $\tilde{U}_1^{(r)}-z\tilde{U}_0^{(r)}.$ In the top part of
fig.\ref{fig11} the estimate of $h(z)$ obtained by Monte
Carlo simulation   is plotted. In the bottom part the
smoothed estimates $\hat{h}(z)$ for $\sigma=0.2$ and
$\beta=5n$ based on a single realization was plotted.

We notice that by the proposed
method we get an improved qualitative information with
respect to that obtained by replicated measures. This is an
important feature for applications where usually only one
data set is measured.

%prodotta da c:\perotti\qualit.m
\begin{figure}
\begin{center}
\hspace{1.cm}{\fbox{\epsfig{file=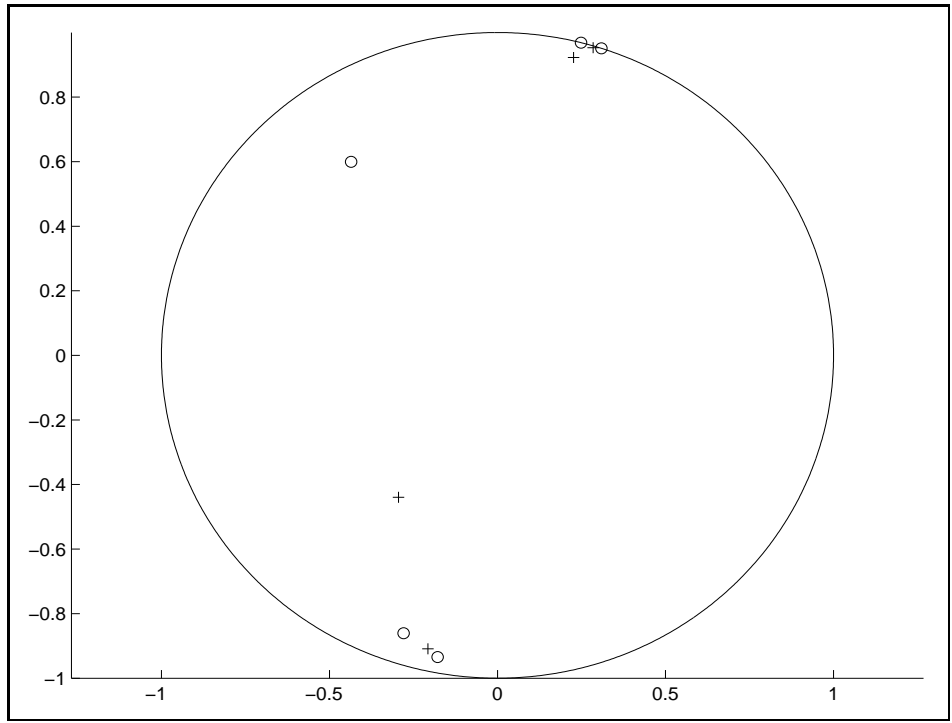,height=9cm,width=12cm}}}
\end{center}
\caption{Poles (circles) and zeros (plus) of $f(z)$ corresponding to  $\underline{\xi}=\left[e^{-0.3-i 2\pi
0.35}, e^{-0.1-i 2\pi  0.3},e^{-0.05-i 2\pi
0.28},e^{-0.0001+i 2\pi 0.2},e^{-0.0001+i 2\pi  0.21}\right],\;\; \underline{c}=\left[ 20,6,3,1,1\right
]$}\label{fig0}
\end{figure}

%prodotta da c:\perotti\figura1.m
\begin{figure}
\begin{center}
%\hspace{1.cm}{\fbox{\epsfig{file=approxzeribb.eps,height=12cm,width=12cm}}}
\hspace{1.cm}{\fbox{\epsfig{file=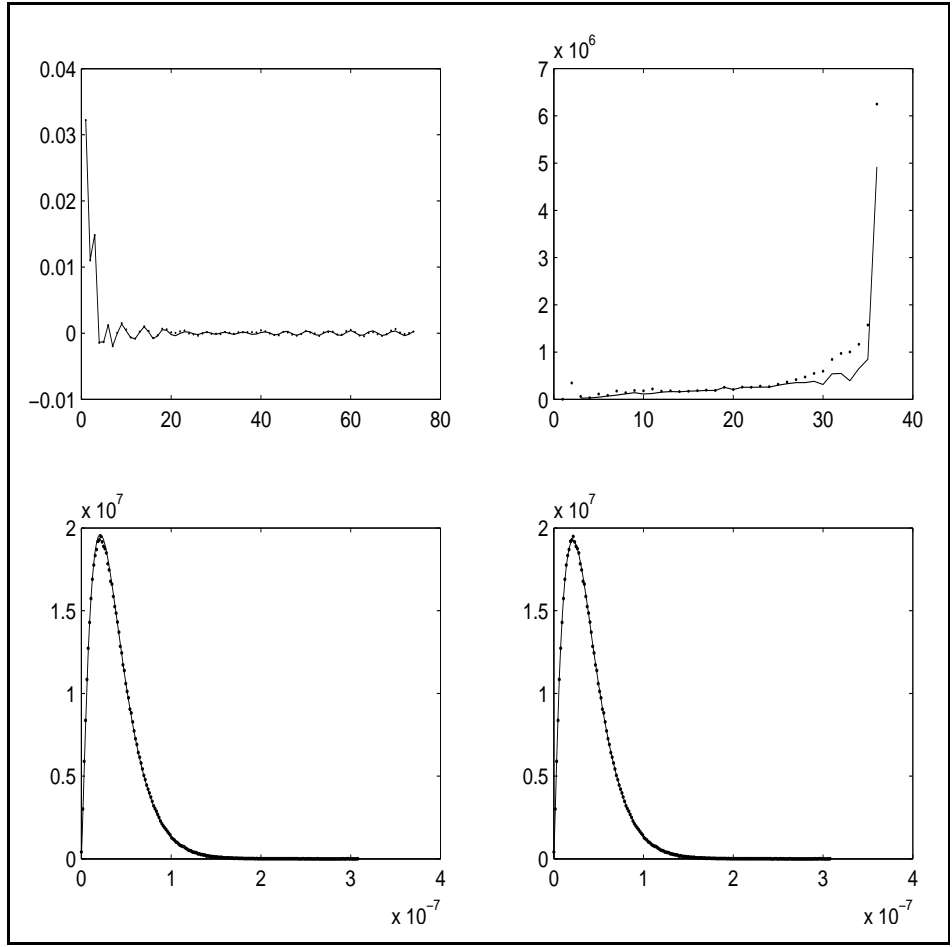,height=12cm,width=12cm}}}
\end{center}
\caption{Top left: real part of the transformed signal (solid) and transformed data (dotted) with $\sigma=0.2$;
 top right:   $L_2$ norm of the difference between the empirical density  of  $\bR_{kk}^2,k=1,\dots,35$  computed by MonteCarlo simulation with $ 10^6$ samples and its approximation obtained by truncating the series expansion of the density after the first term (dotted) and after the first $10$ terms (solid);
  bottom left:  density of $\bR_{kk}^2,\;k=35,$ approximated by the first term of its series expansion (solid), empirical density (dotted);
  bottom right: density of $\bR_{kk}^2,\;k=35,$ approximated by the first $10$ terms of its series expansion (solid), empirical density (dotted).}
 \label{fig0}
\end{figure}

%prodotta da c:\perotti\figura2.m e da figura2biszeri
\begin{figure}
\centering{
\includegraphics[height=10cm,width=15cm]{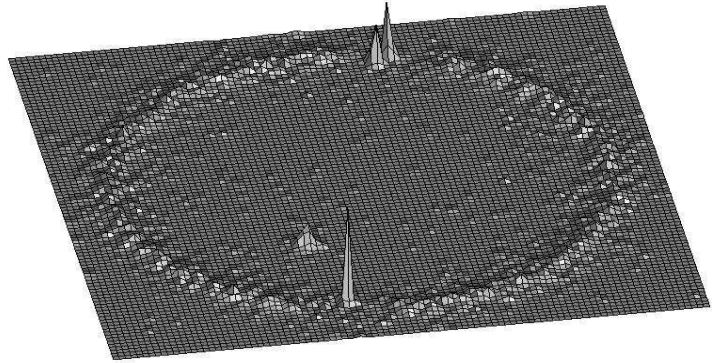}}
\includegraphics[height=10cm,width=15cm]{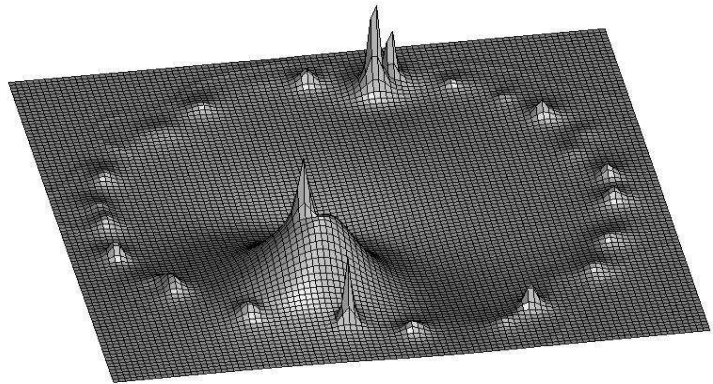}
 \caption{Top: Monte Carlo estimate of the condensed density when $\sigma=0.2$ based on $100$ samples; bottom:
estimate of the condensed density by the closed form approximation
with $\beta=14$ based on one sample.} \label{fig11}
\end{figure}
\end{document}